%% file: ndtcs.tex
\documentclass[aps,pre,twocolumn,groupedaddress,reprint]{revtex4-1}

\usepackage{graphicx}
\usepackage{amsmath}
\usepackage{amsfonts}
\usepackage{amssymb}

\newcommand{\sign}[1] {\mathrm{sign}\left(#1 \right)}
\newcommand{\erfc}[1] {\mathrm{Erfc}\left(#1 \right)}
\newcommand{\erf}[1] {\mathrm{Erf}\left(#1 \right)}
\renewcommand{\Im}[1]{\mathrm{Im}\left(#1 \right)}
\renewcommand{\Re}[1]{\mathrm{Re}\left(#1 \right)}
\newcommand{\mexp}[2][]{\mathbb E^{#1}\! \left[#2 \right]}
\newcommand{\op}[1]{\boldsymbol{\mathrm{#1}}}
\renewcommand{\vec}[1]{\boldsymbol{#1}}

\newtheorem{mythm}{Theorem}

\begin{document}


\title{On energy dissipation in a friction-controlled slide of a body excited by random motions of the foundation}


\author{Sergey Berezin}
\author{Oleg Zayats}
\affiliation{Applied mathematics department, Peter the Great Saint-Petersburg Polytechnic University - Russia, 195251, St.Petersburg, Polytechnicheskaya, 29}


\date{\today}

\begin{abstract}
  We study a friction controlled slide of a body excited by random motions of the foundation it is placed on. Specifically, we are interested in quantities such as displacement, traveled distance, and energy loss due to friction. Assuming that the random excitation is switched off at some time, it is shown that the problem can be treated in an analytic, explicit, manner. Particularly, we derive formulas for the moments of the displacement and distance, and also for the average energy loss. To accomplish that we use the Pugachev--Sveshnikov equation for the characteristic function of a continuous random process given by a system of SDEs. This equation is solved by reduction to a parametric Riemann boundary value problem of complex analysis.
\end{abstract}


\maketitle


\section{\label{sec:intro}Introduction}
In the present paper, we address a problem that concerns dynamical behavior of a solid body sliding with dry friction over the surface of a foundation that moves randomly. This problem is known to be extremely difficult yet very important for various physical applications. 

Interest in such a topic has been increasingly growing over the last few decades. In the 60s, it was engineers who were first to start studying the influence of random fluctuations on the dry friction phenomenon in order to model the behavior of buildings during earthquakes~\cite{Caughey1961, Crandall1974}. More recently, physicists became involved and began to study similar problems but from a ``more microscopic'' perspective. Their studies relate to nanofrictional systems~\cite{Riedo2004}, particles separation~\cite{Eglin2006}, ratchets~\cite{Buguin2006, Fleishman2007, Gnoli2013, Sarrachino2013}, granular motors~\cite{Talbot2011, Gnoli2013a}, and dynamics of droplets on a moving surface~\cite{Daniel2005, Mettu2010, Goohpattader2010}. The quality all these studies share is that, in a way, they are all connected with forces similar to dry friction.

The typical mechanical problem statement we deal with is given in~\cite{Gennes2005}. The author of that paper studies the motion of an object, a solid body, sliding with friction over the surface of a horizontal uniformly rough foundation, which is vibrating laterally being subjected to an external Gaussian white noise excitation. The author uses the simplest description possible for dry friction, namely, he takes into account only the kinetic friction, ignoring the static one. In this case, if we take the lower body as a frame of reference, relative velocity of the sliding body will satisfy a simple stochastic differential (Langevin) equation with~"$\sign{x}$." The latter means the resistant force has a constant absolute value, but its direction is always opposite to the one of the velocity. This equation was a starting point for numerous mathematical studies on the topic, in turn leading to solution of different practical problems.

So far, authors of physical papers have restricted themselves to calculating only the velocity characteristics of the object. The main approach they used is the Fokker--Planck equation~\cite{Touchette2010}, the path integral~\cite{Baule2010}, or the weak-noise limit~\cite{Chen2013}. At the same time, engineers and the specialists in related fields have always wanted to have a more detailed description of the problem. For instance in the mid 70s, S.~H.~Crandall proposed to study the displacement (See~\cite{Crandall1974}). Soon, the corresponding problem was named after Crandall, but neither he nor his co-authors found the exact solution: they had to be satisfied with the approximate one given by statistical linearization. 

It is worth noticing there exist further mechanical studies of classic Crandall's problem and of more general ones. For instance, some authors improved the statistical linearization technique~\cite{Ahmadi1983}, others changed the type of the excitation used~\cite{Constantinou1984} as well as studied the case of two sliding bodies~\cite{Younis1984}. However, all of these works are, at the end of the day, based on the concept of linearization and give an approximate solution only. We should also mention one of more-recent works~\cite{Chen2014}, in which the backward Kolmogorov equation approach is used. The author finds the exact solution for the displacement in terms of the Laplace transform as well as for some other integral functionals. Later on these Laplace images are used to obtain the long time asymptotics.

In the present paper we suggest an alternative method based on the Pugachev--Sveshnikov equation. This equation is not wide known, yet it is very effective for the problems considered. Unlike the Fokker--Planck equation, the Pugachev--Sveshnikov equation describes the behavior of the random process in terms of the characteristic function. This equation allows us to get an exact expression not only for the object's velocity, but also for the displacement and for the total distance traveled. This distance is proportional to the amount of the energy dissipated during the slide, and as it turns out there is a certain natural conservation law for it.

The Pugachev--Sveshnikov equation method is thoroughly described in~\cite{Zayats2013, Zayats2013a}; also some preliminary study of the present topic can be found in~\cite{Zayats2007}.

\section{\label{sec:probl_st}Problem statement}
We consider a rigid body of mass~$m$ placed on a massive foundation. This foundation is subjected to the random Gaussian white noise excitation~$\xi(t)$ of intensity~$h$, switched off at time~$t_0$; the corresponding covariance function is~$K_\xi(t_1, t_2) = \delta(t_2 - t_1)$. That excitation causes body to move with velocity~$V(t)$ relative to the foundation.

\begin{figure}[h!]
  \includegraphics[width = 8.5cm]{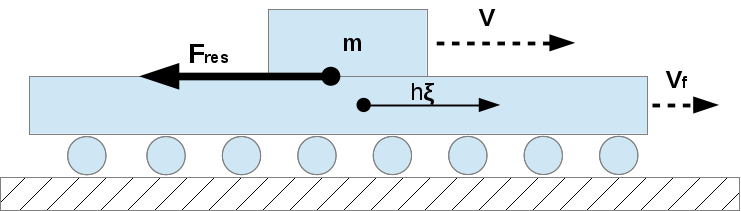}
  \caption{Crandall's problem}
  \label{fig1}
\end{figure}

The resistant force~$F_{res}(t)$ between the foundation and the body is assumed to obey Coulomb's friction law, namely, $F_{res} (t)= -\mu m g \, \sign{V(t)}$, where~$\mu$ is the coefficient of dry (kinetic) friction between the two surfaces in contact (see Fig.~\ref{fig1}), and~$g$ is the acceleration of gravity.

The body's velocity~$V(t)$, displacement~$U(t)$, and the distance~$S(t)$ body has traveled satisfy equations that read as follows:
\begin{equation}
  \begin{aligned}
    \label{eq1}
    &\dot V(t) = -\mu g \, \sign{V(t)} + \eta(t_0 - t) \, h \xi(t),\\
    &\dot U(t) = V(t), \quad \dot S(t) = |V(t)|,
  \end{aligned}
\end{equation}
where~$\eta(t)$ is the Heaviside step function that is in charge of switching excitation off. The first component~$V(t)$ of the process~\eqref{eq1} is the so-called Caughey--Dienes process~\cite{Caughey1961}, the second component~$U(t)$ is the Crandall process~\cite{Crandall1974}, and the third component~$S(t)$ is a new one, that has not been thoroughly studied to date.

Suppose that the system had been at rest before the moment~$t=0$ when excitation suddenly occurred, and this is when we started measuring~$U(t)$ and~$S(t)$. In such a situation initial conditions become homogeneous~$V(0)=U(0)=S(0)=0$.

System~\eqref{eq1} includes two nonlinearities
\begin{equation}
  \label{eq1a}
  \Psi_1(V) = \sign{V}, \quad \Psi_2(V) = |V|,
\end{equation}
both of which are piecewise linear functions that have two domains of linearity. Therefore, this system can be treated with methods from~\cite{Zayats2013}.

\section{\label{sec:av_energ}Average energy conservation law}
Before diving into direct solution of the problem, we would like to present some alongside result on the average energy conservation.

From now on, we denote specific kinetic energy (per unit mass) of the body by~$Q(t)$ and specific energy dissipated due to friction by~$R(t)$:
\begin{equation}
  \label{eq4}
  \begin{aligned}
    &Q(t) = \frac{1}{2} V^2(t),\\
    &R(t) = \mu g \int \limits_{0}^{t} V(s) \, \sign{V(s)}\, ds = \mu g \, S(t).
  \end{aligned}
\end{equation}
As long as external excitation is switched off, the body stops at some random moment~$t_s > t_0$ because of friction, and thus~$V(t_s)=0$. Let us write the It\^o's formula for~$Q(t)$, using~\eqref{eq1}, in the following form:
\begin{equation}
  \begin{aligned}
    \label{eq4a}
    dQ(t) &= \Big(- \mu g |V(t)| + \frac{1}{2} \eta(t_0-t) h^2 +\\
    &+\eta(t_0-t) h V(t)\, \xi(t) \Big)\, dt.    
  \end{aligned}
\end{equation}
Having taken the mathematical expectation of~\eqref{eq4a}, we will arrive at
\begin{equation}
  \label{eq4b}
  d \bar q(t) = - d \bar r(t) + \frac{1}{2} \eta(t_0-t) h^2 \, dt,
\end{equation}
where $\bar q(t) = \mexp{Q(t)}$, $\bar r(t) = \mexp{R(t)}$, and $\bar{r}_s = \mexp{R(t_s)}$.

Now, if we consider two time intervals~$(0, t_0)$ and~$(t_0, t_s)$, and integrate~\eqref{eq4b} over each of them, taking into account initial condition, we will get the average energy conservation law
\begin{align}
  \label{eq5}
  \bar q(t_0)  &= -\bar w(t_0) + \frac{1}{2} h^2 t_0,\\
  \label{eq6}
  -\bar q(t_0) &= -(\bar w_s-\bar w(t_0)).
\end{align}
The specific kinetic energy of the foundation which moves with velocity~$V_f(t)$ is given by
\begin{equation}
  K_f(t) = \frac{1}{2} V_f^2(t).  
\end{equation}
Since the excitation in~\eqref{eq1} is a Gaussian white noise of intensity~$h$, the process $V_f(t)$ will be the Wiener one. Therefore, the average kinetic energy of the foundation has the form
\begin{equation}
  \bar k_f(t) = \frac{1}{2} h^2 t. 
\end{equation}
That allows us to rewrite~\eqref{eq5} as
\begin{align}
  \label{eq5a}
  \bar q(t_0) &= -\bar w(t_0) + \bar k_f(t_0),\\
  \label{eq6a}
  -\bar q(t_0) &= -(\bar w_s-\bar w(t_0)).
\end{align}
Summing up expressions~\eqref{eq5a} and~\eqref{eq6a}, we will finally get to the formula
\begin{equation}
  \label{eq7}
  \bar w_s = \bar k_f(t_0),
\end{equation}
which says that the average specific energy spent on the friction during the motion equals the average specific kinetic energy of the foundation at the time the excitation is switched off. This result give us an opportunity to calculate the average energy due to friction at the final moment~$t_s$. To investigate transient behavior we are going to use the Pugachev--Sveshnikov equation approach. The latter let us find the moments of velocity~$V(t)$, displacement~$U(t)$, distance~$S(t)$, and, therefore, these of dissipated energy~$R(t)$.

\section{\label{sec:PS_eq}Pugachev--Sveshnikov equation formalism}
As we pointed out earlier, system~\eqref{eq1} is a piecewise linear one and has two domains of linearity. This allows us to use the Pugachev--Sveshnikov equation approach from~\cite{Zayats2013}.

To simplify the following calculations, we scale equations~\eqref{eq1} to dimensionless form using the transform
\begin{equation}
  \label{eq8}
  \begin{aligned}
    &U_1 = \frac{\mu g}{h^2} V,\quad U_2 = \frac{\mu^3 g^3}{2 h^4} U,\\
    &U_3 = \frac{\mu^3 g^3}{2 h^4}S,\quad \tau = \frac{\mu^2 g^2}{2 h^2}t.
  \end{aligned}
\end{equation}
Now we can rewrite~\eqref{eq1} as
\begin{equation}
  \label{eq3}
  \begin{aligned}
    &\dot U_1(\tau) = -2\, \sign{U_1(\tau)} + \eta(\tau_0 - \tau)\, \sqrt{2} \tilde{\xi}(\tau),\\
    &\dot U_2(\tau) =  U_1(\tau),\quad \dot U_3(\tau) = |U_1(\tau)|.
  \end{aligned}  
\end{equation}
The process~$\tilde{\xi}(\tau)$ is another Gaussian white noise process given by
\begin{equation}
  \tilde{\xi}(\tau) = \frac{\sqrt{2} h}{\mu g}\, \xi{\left(\frac{2 h^2}{\mu^2 g^2} \tau \right)}.
\end{equation}

In what follows, we review the main steps of the Pugachev--Sveshnikov equation approach in order to find the characteristics of~$(U_1(\tau),U_2(\tau),U_3(\tau))$. First, assuming~$t \le t_0$, the singular integral-differential equation for the characteristic function~$E(z_1,z_2,z_2; \tau)$ of the Markov process~$(U_1(\tau), U_2(\tau), U_3(\tau))$ will take the form (see the appendix)
\begin{equation}
  \label{eq10}
  \begin{aligned}
    &\frac{\partial E}{\partial \tau} - z_2 \frac{\partial E}{\partial z_1} + z_1^2 E + \frac{2 z_1}{\pi} \mathrm{v.p.} \int \limits_{-\infty}^{+\infty} \frac{E|_{z_1=s}}{s - z_1} ds +\\
    &+\frac{i z_3}{\pi} \frac{\partial}{\partial z_1} \left[ \mathrm{v.p.} \int \limits_{-\infty}^{+\infty} \frac{E|_{z_1=s}}{s - z_1} ds\right]=0
  \end{aligned}
\end{equation}
with~$E|_{\tau=0} = 1$. It is important to underline that the integral in~\eqref{eq10} is improper, the principle value one.

To solve equation~\eqref{eq10}, let us introduce the Cauchy type integral for all~$\zeta \in \mathbb{C}$ such that~$\Im{\zeta} \ne 0$:
\begin{equation}
  \label{eq10a}
  F(\zeta; z_2, z_3, \tau) = \frac{1}{2 \pi i} \int \limits_{-\infty}^{+\infty} \frac{E(s, z_2, z_3; \tau)}{s - \zeta} ds,
\end{equation}
which is known to be a complex piecewise analytic function. For~$z_1 \in \mathbb{R}$, we denote the limit values of~$F$ as~$\zeta \to z_1 \pm 0 i$ by~$F^\pm(z_1; z_2, z_3, \tau)$. The latter functions satisfy the well-known Sokhotski--Plemelj formulas~\cite{Gakhov1990}
\begin{equation}
  \label{eq10b}
  E = F^{+} - F^{-},\quad \frac{1}{\pi i} \mathrm{v.p.} \int \limits_{-\infty}^{+\infty} \frac{E|_{z_1=s}}{s - z_1} ds  = F^{+} + F^{-}.
\end{equation}

We insert~\eqref{eq10b} into~\eqref{eq10} and rearrange the terms in the following way:
\begin{equation}
  \label{eq10bb}
  \begin{aligned}
    &\frac{\partial F^+}{\partial \tau} - (z_2 + z_3)\frac{\partial F^+}{\partial z_1} + z_1 (z_1 + 2 i) F^+ =\\
    &\frac{\partial F^-}{\partial \tau} - (z_2 - z_3)\frac{\partial F^-}{\partial z_1} + z_1 (z_1 - 2 i) F^-.
  \end{aligned}
\end{equation}
Then notice, that the left-hand side of~\eqref{eq10bb} can be analytically continued to the upper half plane with respect to~$z_1$, the right-hand side to the lower half plane, and for all~$z_1 \in \mathbb{C}$ such that~$\Im{z_1}=0$ the both sides match. This means that the left- and right-hand sides of~\eqref{eq10bb} are elements of the same entire analytic function. Assuming that~$F^\pm = O(1/z_1)$ as~$z_1 \to \infty$ such that~$\Im{z_1}\gtrless 0$, we have that this entire function grows no faster than linearly. Thus, due to Liouville's theorem from complex analysis, this entire function is, in fact, linear, and we get the so-called master equation
\begin{equation}
  \label{eq10c}
  \begin{aligned}
    &\frac{\partial F^\pm}{\partial \tau} - (z_2 \pm z_3)\frac{\partial F^\pm}{\partial z_1} + z_1 (z_1 \pm 2 i) F^\pm =\\
    &=G_0 + z_1 G_1,\quad F^\pm|_{\tau=0} = \pm \frac{1}{2}.
  \end{aligned}
\end{equation}
Here, $G_0(z_2, z_3; \tau)$ and~$G_1(z_2, z_3; \tau)$ are some new intermediary complex functions we need to find. Notice that, as follows from~\cite{Berezin2016}, the characteristic function~$E(z_1,z_2,z_3; \tau)$ is entire with respect to~$(z_1,z_2,z_3)$. Similar argument shows that~$F^\pm(z_1; z_2, z_3; \tau)$ are entire with respect to~$(z_2, z_3)$ as well, and, therefore, the same goes for~$G_0$ and~$G_1$ due to~\eqref{eq10c}.

In the next section we show how to find~$G_0$ and~$G_1$. Now, we just mention that after these functions are found, we can solve equations~\eqref{eq10c}, and use~\eqref{eq10b} to find the desired characteristic function~$E$.

\section{\label{sec:moments}Method of moments}
For the sake of simplicity and to fulfill practical needs, we find moments of~$(U_1(\tau), U_2(\tau), U_3(\tau))$ only, although, notice that the distribution itself can be found in the similar manner (See~\cite{Zayats2013}). Now, we review the basic steps from~\cite{Zayats2013}. First, let us start by applying the Laplace transform to~\eqref{eq10c} with respect to~$\tau$. The corresponding equations read
\begin{equation}
  \label{eq10d}
  - (z_2 \pm z_3)\frac{\partial \tilde F^\pm}{\partial z_1} + (z_1^2 \pm 2 i z_1 + p) \tilde F^\pm = \tilde G_0 + z_1 \tilde G_1,  
\end{equation}
where~$\tilde G_{0}(z_2, z_3; p), \tilde G_1(z_2, z_3; p)$, and $\tilde F^\pm(z_1; z_2, z_3, p)$ are the Laplace transforms of~$G_{0}(z_2, z_3; \tau), G_1(z_2, z_3; \tau)$, and~$F^\pm(z_1; z_2, z_3, \tau)$, correspondingly. These transforms exist for~$\Re{p} > 0$, and, preserve all the analytic properties of the originals at least for the small neighborhood of~$(z_2,z_3)=(0,0)$. The corresponding formal series decomposition reads
\begin{align}
  \label{eq10e}
  &\tilde F^\pm(z_1; z_2, z_3, p) = \sum \limits_{\alpha=0}^\infty \sum \limits_{\beta=0}^\infty \tilde F^\pm_{\alpha \beta} (z_1; p) \frac{z_2^\alpha z_3^\beta}{\alpha ! \beta !}\\
  \label{eq10f}
  &\tilde G_j(z_2, z_3, p) = \sum \limits_{\alpha=0}^\infty \sum \limits_{\beta=0}^\infty \tilde G_{j,\alpha \beta} (p) \frac{z_2^\alpha z_3^\beta}{\alpha ! \beta !},
\end{align}
where~$\tilde F^+_{\alpha \beta} (z_1; p)$ and~$\tilde F^-_{\alpha \beta} (z_1; p)$ are analytic in the upper- and lower-half planes respectively.

Next, inserting~\eqref{eq10e} and~\eqref{eq10f} into~\eqref{eq10d} and comparing coefficients of~$z_2^k z_3^l$ on both sides, we will have the infinite system of equations
\begin{equation}
  \label{eq10g}
  \begin{aligned}
    &(z_1^2 \pm 2 i z_1 + p) \tilde F^\pm_{kl} = \tilde G_{0,kl} + z_1 \tilde G_{1,kl} \pm \frac{1}{2} \delta_{0k} \delta_{0l} +\\
    &k \frac{\partial \tilde F^\pm_{(k-1)l}}{\partial z_1} \pm l \frac{\partial \tilde F^\pm_{k(l-1)}}{\partial z_1}, \quad k,l \ge 0,
  \end{aligned}
\end{equation}
where~$\delta$ is the Kronecker delta, and~$\tilde F^\pm_{kl} =0$ if either~$k<0$ or~$l<0$.

To get mixed moments of~$(U_1(\tau),U_2(\tau),U_3(\tau))$ we introduce the functions
\begin{equation}
  \label{eq10h}
  m_{jkl} (\tau) = \mexp{U_1^j(\tau) U_2^k(\tau) U_3^l(\tau)}.
\end{equation}
Formulas~\eqref{eq10b} and~\eqref{eq10e} lead us to the expression
\begin{equation}
  \label{eq10i}
  \tilde m_{jkl} (p) =\frac{1}{i^{j+k+l}} \frac{\partial^j}{\partial z_1^j} \left( \tilde F^+_{kl}(0; p) - \tilde F^-_{kl}(0; p) \right). 
\end{equation}

We will calculate~$\tilde F^\pm_{kl}$ successively by grouping them with respect to~$k$ and~$l$, so that~$k+l=0,1,2,...$. There is the one term for~$k+l=0$, the term~$\tilde F^\pm_{00}$,  which can be easily found from~\eqref{eq10g}:
\begin{equation}
  \label{eq10j}
  \tilde F^\pm_{00} = \frac{\tilde G_{0, 00} + z_1 \tilde G_{1, 00} \pm \frac{1}{2}}{z_1^2 \pm 2 i z_1 + p}.
\end{equation}
These rational functions have simple poles in the complex plane, and in order to turn them into analytic functions for~$\Im{z_1} \gtrless 0$, respectively, we need to vanish enumerators of the fractions at~$\pm i \mu$, where~$\mu=\sqrt{p+1}-1$. That gives us a system of linear equations to find~$\tilde G_{0,00}$ and~$\tilde G_{1,00}$:
\begin{align}
  \label{eq10k}
  \tilde G_{0, 00} + i \mu \tilde G_{1, 00} + \frac{1}{2} = 0,\\
  \label{eq10kk}
  \tilde G_{0, 00} - i \mu \tilde G_{1, 00} - \frac{1}{2} = 0.
\end{align}
And finally we get
\begin{equation}
  \label{eq10l}
  \tilde F^\pm_{00} = \frac{i}{2 \mu (z_1 \pm i (\sqrt{p+1} + 1))}.
\end{equation}
For~$k+l > 0$ we will have the following recurrent formulas
\begin{equation}
  \label{eq10m}
  \tilde F^\pm_{kl} = \frac{\tilde G_{0,kl} + z_1 \tilde G_{1,kl}+ k \frac{\partial \tilde F^\pm_{(k-1)l}}{\partial z_1} \pm l \frac{\partial \tilde F^\pm_{k(l-1)}}{\partial z_1}}{z_1^2 \pm 2 i z_1 + p}.
\end{equation}

In the similar way as earlier, to make these complex valued functions analytic, we need to vanish enumerators at~$i \mu^\pm$. That gives us two other equations
\begin{equation}
  \label{eq10n}
  \begin{aligned}
    &\tilde G_{0,kl} \pm i \mu \tilde G_{1,kl}+ k \frac{\partial \tilde F^\pm_{(k-1)l}(\pm i \mu; p)}{\partial z_1} \pm\\
    &\pm l \frac{\partial \tilde F^\pm_{k(l-1)}(\pm i \mu; p)}{\partial z_1} = 0,
  \end{aligned}
\end{equation}
where we can take either~$+$ or~$-$ for all of~$\pm$ simultaneously.

Using the procedure  described, we will easily get Laplace transforms of the first moments:
\begin{align}
  \label{eq10o}
  &\tilde m_{100} = \tilde m_{010} = 0, \quad \tilde m_{001}=\frac{1}{p^2 (\sqrt{p+1}+1)},\\
  &\tilde m_{200}=\frac{2}{p (\sqrt{p+1}+1)^2},\\
  \label{eq10p}
  &\tilde m_{020}=\frac{4\sqrt{p+1} + 1}{p^2 (p+1) (\sqrt{p+1}+1)^3},\\
  \label{eq10pp}
  &\tilde m_{002}=\frac{9 p - 4\sqrt{p+1} + 8}{2 p^3 (p+1) (\sqrt{p+1}+1)^2}.
\end{align}
All the rest moments can be found in the similar manner.

\section{\label{sec:result}Results}
Applying the inverse Laplace transform to~\eqref{eq10o}~--~\eqref{eq10pp} we will have
\begin{equation}
  \label{eq10q}
  \bar u_1(\tau) = 0, \quad \bar u_2(\tau) = 0,
\end{equation}
\begin{widetext}
  \begin{equation}
    \label{eq10r}
    \bar u_3(\tau) = \frac{1}{8} \left[2(1+2 \tau)\sqrt{\frac{\tau}{\pi}} e^{-\tau} - 4 \tau^2 \erfc{\sqrt{\tau}}+(4 \tau -1) \erf{\sqrt{\tau}} \right],
  \end{equation}
  \begin{equation}
    \label{eq10s}
    \sigma_{U_1}^2(\tau) = \frac{1}{2} \left[ \erf{\sqrt{\tau}} + 4 \tau (\tau + 1) \erfc{\sqrt{\tau}} -2(1+2 \tau) \sqrt{\frac{\tau}{\pi}} e^{-\tau} \right],
  \end{equation}
  \begin{equation}
    \label{eq10t}
    \sigma_{U_2}^2(\tau) = \frac{5}{8} \tau - \frac{27}{32} + e^{-\tau}\left[1-\sqrt{\frac{\tau}{\pi}} \left( \frac{1}{2}\tau^3 +\frac{7}{12}\tau^2-\frac{13}{24} \tau + \frac{5}{16} \right) \right] + \erfc{\sqrt{\tau}} \left( \frac{1}{2} \tau^4+\frac{5}{6} \tau^3-\frac{1}{2}\tau^2 + \frac{3}{8}\tau -\frac{5}{32} \right),
  \end{equation}
  \begin{equation}
    \label{eq10u}
    \sigma_{U_3}^2(\tau) = \frac{1}{4}\tau^2 + \frac{1}{8}\tau - \frac{11}{32} + e^{-\tau}\left[\frac{1}{2} +\tau \sqrt{\frac{\tau}{\pi}} \left(\frac{13}{24} -\frac{7}{12}\tau - \frac{1}{2}\tau^2\right)\right] + \erfc{\sqrt{\tau}} \left( \frac{1}{2} \tau^4 + \frac{5}{6} \tau^3 -\frac{1}{2}\tau^2 + \frac{3}{8}\tau -\frac{5}{32} \right)  -\bar u_3^2    
  \end{equation}
\end{widetext}
for $\tau \le \tau_0 = \frac{\mu^2 g^2}{2 h^2}t_0$, where $\erf{x} = \frac{2}{\sqrt{\pi}} \int \limits_0^{x} e^{-s^2} ds$ and~$\erfc{x} = 1 - \erf{x}$.

It is also worth noticing that the expression for~$\tilde m_{020}$ in~\eqref{eq10p} perfectly matches that in formula~(73) of paper~\cite{Chen2014}, having in mind the proper scaling. The same goes for the long time asymptotics of~$\sigma^2_{U_2}(\tau)$ in~\eqref{eq10t} (compare with formula~(76) in~\cite{Chen2014}).

After the external excitation is turned off, the body will remain sliding until the moment~$t_s$, or $\tau_s=\frac{\mu^2 g^2}{2 h^2}t_s$ in the dimensionalless formulation. The expression~\eqref{eq10l} gives us a possibility to derive the probability density function of~$\tau_s$. Inverting~\eqref{eq10l} by Fourier with respect to~$z_1$ and by Laplace with respect to~$p$, one can get
\begin{equation}
  \begin{aligned}
    f_{\tau_s}(y, \tau) &=  \frac{1}{2} \left( e^{-2 |y|} \erfc{\frac{|y|-2 \tau}{2 \sqrt{\tau}}}+ \right.\\
    &\left. + \frac{1}{\sqrt{\pi \tau}} e^{-\frac{(|y|+2\tau)^2}{4\tau}} \right).
  \end{aligned}
\end{equation}

Since for~$\tau > \tau_0$ the equations~\eqref{eq3} will have no stochastic component, it is easy to solve them. Then, we will have
\begin{align}
  \label{eq13a}
  U_1(\tau) &= U_1(\tau_0) - 2 \sign{U_1(\tau_0)}\; (\tau-\tau_0),\\
  \label{eq13b}
  U_2(\tau) &= U_2(\tau_0) + U_1(\tau_0) (\tau - \tau_0) -\\
            &-\sign{U_1(\tau_0)}\; (\tau -\tau_0)^2,\\
  \label{eq13c}
  U_3(\tau) &= U_3(\tau_0) + |U_1(\tau_0)| (\tau-\tau_0) - (\tau-\tau_0)^2,
\end{align}
where~$\tau \le \tau_s$, and at time~$\tau_s$ the body stops.

From~\eqref{eq13a} it follows that~$\tau_s$ is given by the formula
\begin{equation}
  \label{eq14}
  \tau_s = \tau_0 + \frac{1}{2}|U_1(\tau_0)|.
\end{equation}

If we substitute~\eqref{eq14} into expression~\eqref{eq13c}, go back to units with~\eqref{eq8}, and use~\eqref{eq4}, we will have the formula~$\bar r_s = \bar r(t_0) + \bar q(t_0)$. Thus, by use of~\eqref{eq10r} and~\eqref{eq10s} after units conversion~\eqref{eq8}, we will have~$\bar r_s = \frac{h^2 t_0}{2} = \bar k_f(t_0)$, which is in perfect agreement with~\eqref{eq7}.

\begin{figure}[h!]
  \includegraphics[width = 6cm]{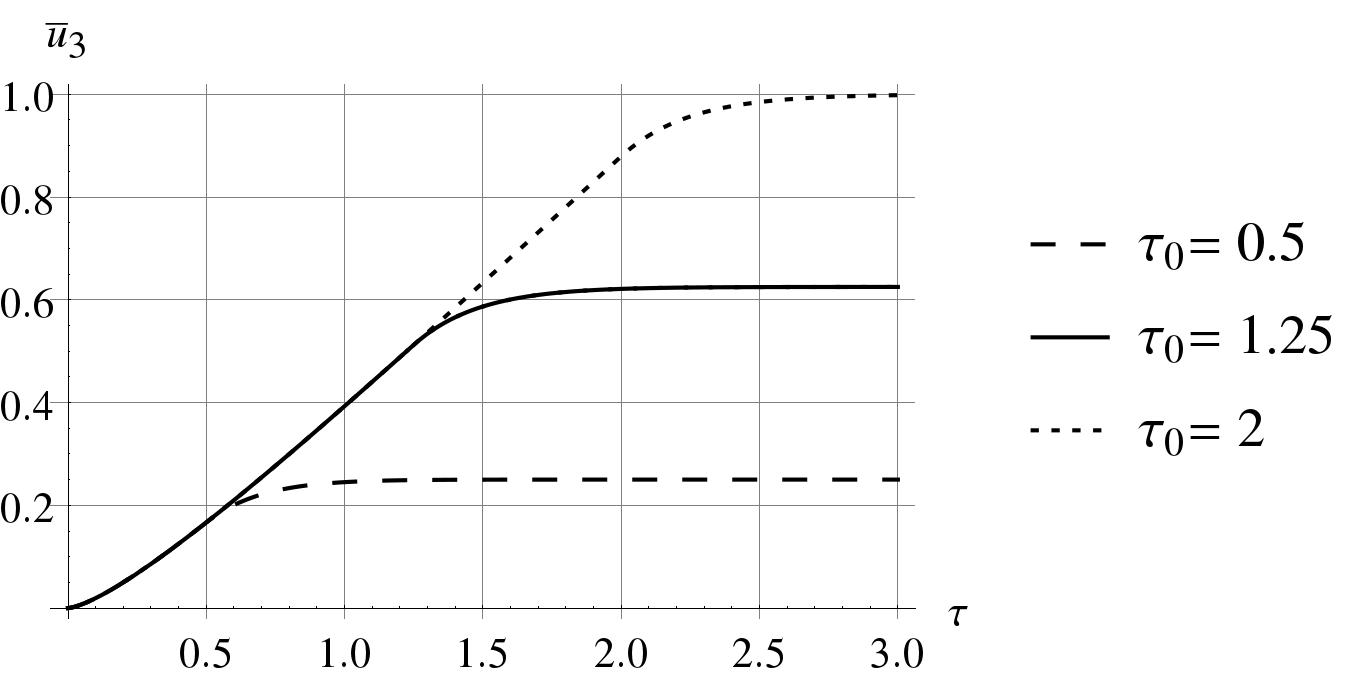}
  \caption{Average scaled traveled distance~$\bar u_3(\tau)$}
  \label{fig2}
\end{figure}

\begin{figure}[h!]
  \includegraphics[width = 6cm]{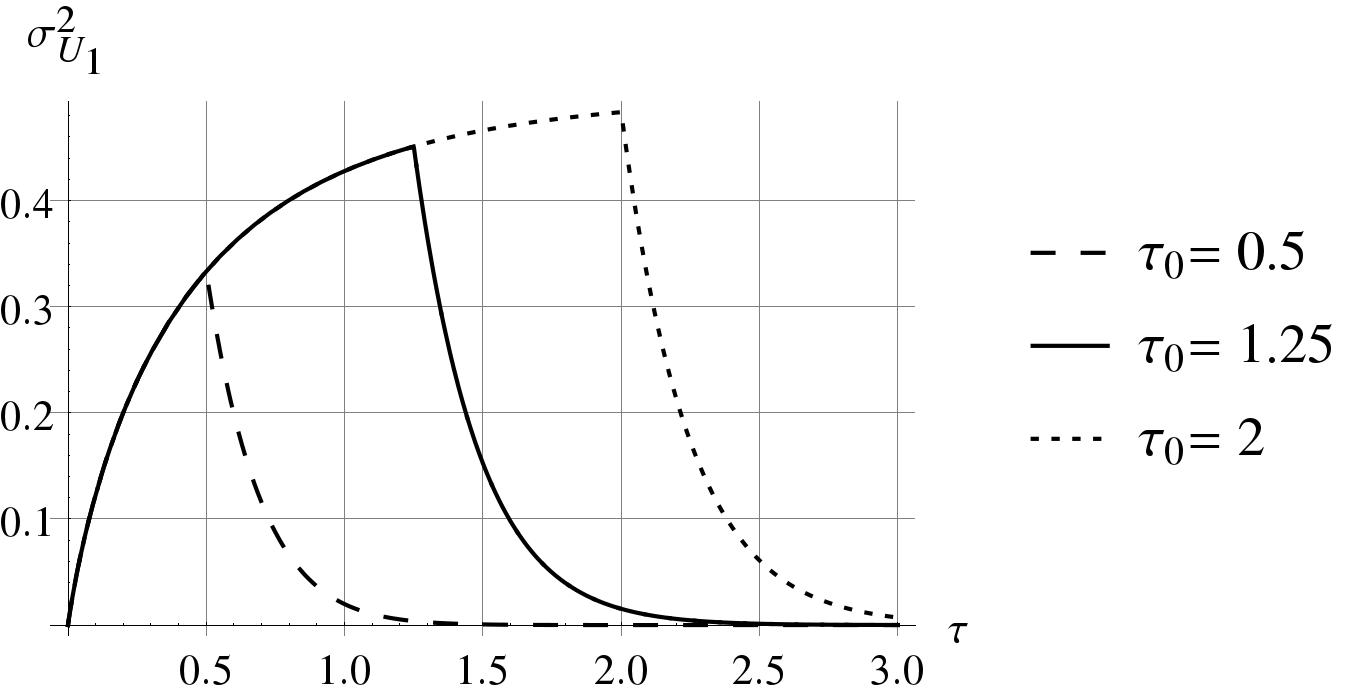}
  \caption{Variance of scaled velocity~$\sigma^2_{U_1}(\tau)$}
  \label{fig3}
\end{figure}

\begin{figure}[h!]
  \includegraphics[width = 6cm]{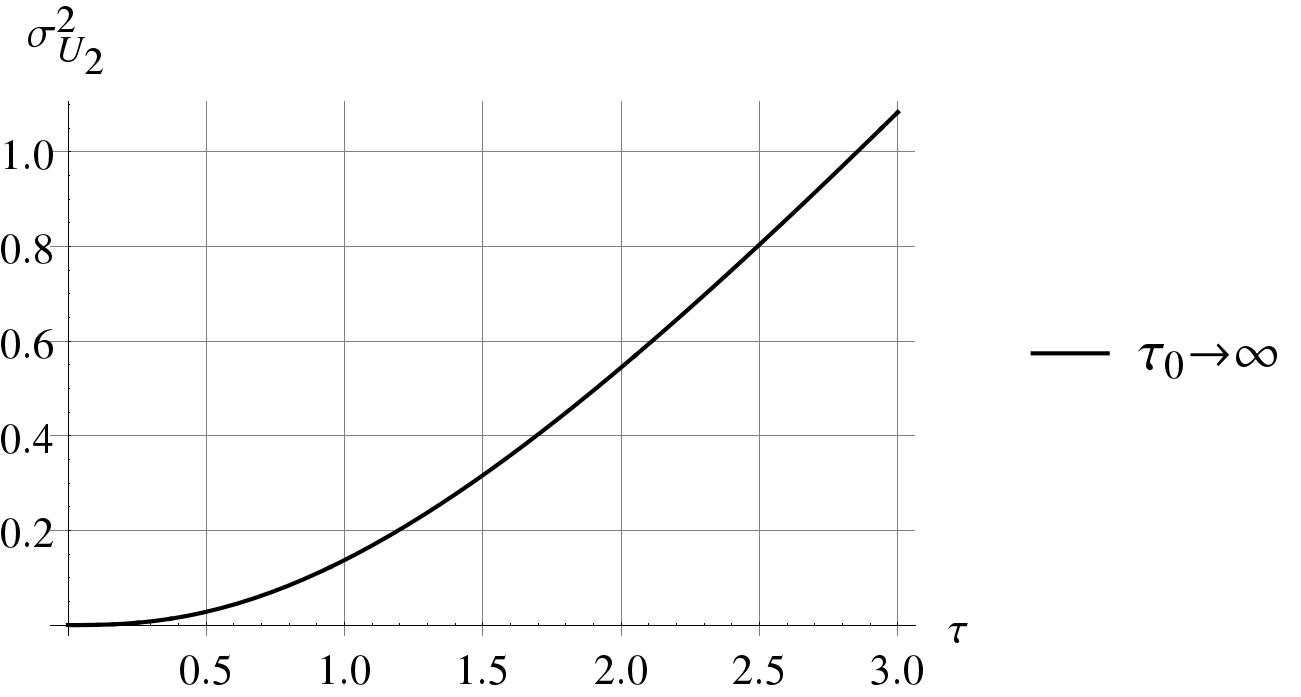}
  \caption{Variance of scaled displacement~$\sigma^2_{U_2}(\tau)$}
  \label{fig4}
\end{figure}

\begin{figure}[h!]
  \includegraphics[width = 6cm]{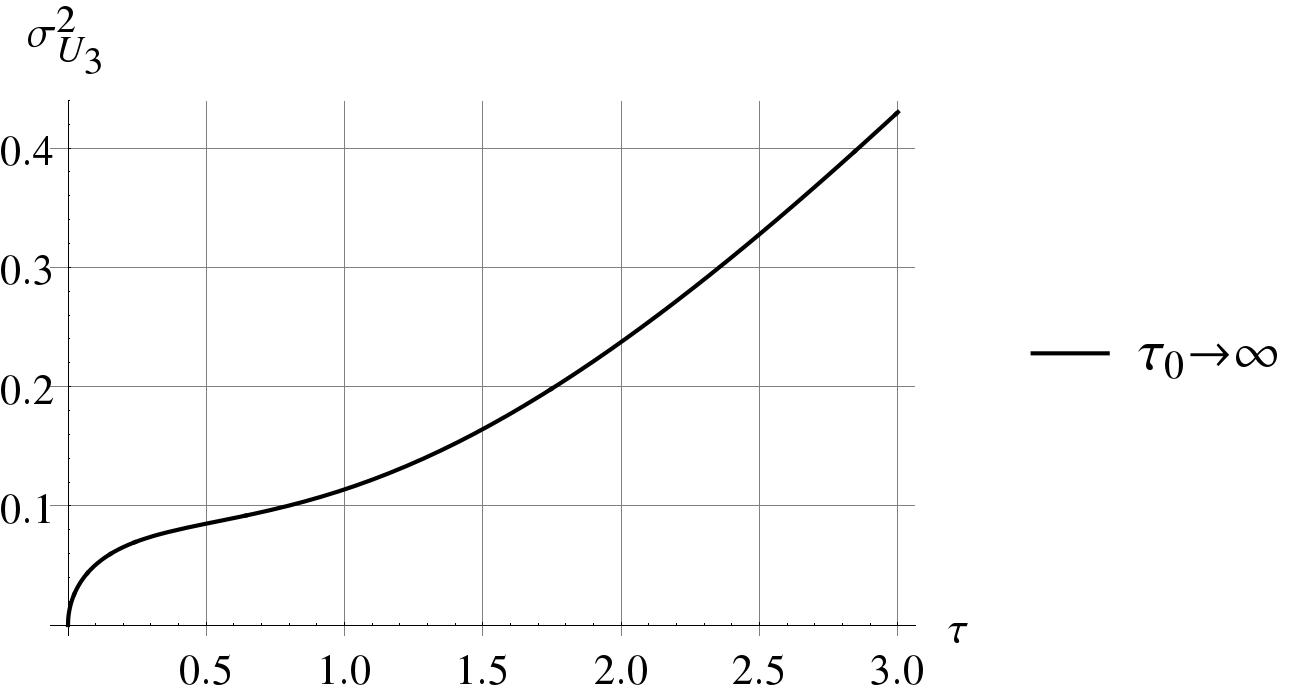}
  \caption{Variance of scaled traveled distance~$\sigma^2_{U_3}(\tau)$}
  \label{fig5}
\end{figure}

Finally, the formulas~\eqref{eq10q}--\eqref{eq14} allow us to plot some interesting pictures (see Fig.~\ref{fig2}, \ref{fig3}, \ref{fig4}, and \ref{fig5}).



\begin{acknowledgments}
  The authors would like to thank Alessandro Sarracino, Andrea Gnoli, and Andrea Puglisi, as well as the anonymous referees, for the interest in this work, and for additional papers on the subject they kindly sent to us. There is no doubt that strengthened the overall quality of our article.
\end{acknowledgments}

\appendix*
\section{\label{app:PS_eq}The Pugachev--Sveshnikov equation}
In the appendix for the sake of completeness, we recap main points of the Pugachev--Sveshnikov equation theory used above. Its rigorous consideration can be found in~\cite{Zayats2013, Zayats2013a}. To be systematic, we start from a general vector stochastic differential equation in~$\mathbb R^n$:
\begin{equation}
  \begin{aligned}
    \label{eq:SDE_gen}
    &d \vec{X}(t) = \vec{a}(t, \vec{X}(t)) dt +\op{H}(t, \vec{X}(t))\, d\vec{W}(t),\\
    &\vec{X}(0)=\vec{X}_0,
  \end{aligned}
\end{equation}
where~$\vec{X}(t)$ is a vector process, $\vec{a}(t, \vec{x})$ is the drift vector, $\op{H}(t, \vec{x})$ is the noise intensity matrix, and $\vec{W}(t)$ is the vector of independent Wiener processes. In applications the latter is usually written as
\begin{equation}
  \label{eq:SDE_gen_Lang}
  \dot{\vec{X}}(t) = \vec{a}(t, \vec{X}(t)) +\op{H}(t, \vec{X}(t))\, \vec{\xi}(t),
\end{equation}
where~$\vec{\xi}(t)$ is the standard vector Gaussian white noise.

Under some mild assumptions on the coefficients~$\vec{a}(t, \vec{x})$,~$\op{H}(t, \vec{x})$ and the initial random vector~$\vec{X}_0$ it is known that there exist a unique solution of~\eqref{eq:SDE_gen}, that is a Markov process. It is also well known that the probability density function~$f_X(\vec{x}, t)$ satisfies Fokker--Planck--Kolmogorov equation. In mid 40s it was shown that not only is there an equation for the density, but also for the characteristic function
\begin{equation}
  E(\vec{z}, t) = \mexp{e^{i \vec{z}^* \vec{X}(t)}}, z\in \mathbb{R}^n,
\end{equation}
where the operation~$(\cdot)^*$ is the transpose. Such an equation reads
\begin{widetext}
  \begin{equation}
    \label{eq:Pug_gen}
    \frac{\partial E(\vec{z},t)}{\partial t} = \mexp{\bigg(i \vec{z}^* \vec{a}(t, \vec{X}(t)) - \frac{1}{2} \vec{z}^* \op{B}(t, \vec{X}(t)) \vec{z}\bigg) e^{i \vec{z}^* \vec{X}(t)}}, \quad E(\vec{z},0) = E_0(\vec{z}).
  \end{equation}
\end{widetext}
Here, the function~$E_0(\vec{z})$ is the characteristic function of~$\vec{X}_0$, given by~$E_0(\vec{z})=\mexp{e^{i z^* \vec{X}_0}}$.
The diffusion matrix~$\op{B}(t, \vec{x})$ is defined in the following manner: $\op{B}(t, \vec{x}) = (\op{H}(t, \vec{x})) (\op{H}(t, \vec{x}))^*$.

The above equation easily follows from It\^o's formula for~$e^{i \vec{z}^* \vec{X}(t)}$ and in fact can be modified for more general types of SDEs (e.g., SDEs driven by Levy processes, etc.). Rigorous information on this equation and its applications can be found in monograph~\cite{Pugachev1987a}.

It is worth noticing that equation~\eqref{eq:Pug_gen} is an equation with respect to the characteristic function~$E(\vec{z}, t)$. Indeed, the mathematical expectation in the right hand side of~\eqref{eq:Pug_gen} is completely defined by the one-dimensional distribution of~$\vec{X}(t)$, which is, in its turn, completely defined by the characteristic function~$E(\vec{z},t)$.

Quite for the long time, equation~\eqref{eq:Pug_gen} had been thought of as difficult to solve explicitly, and consequently was abandoned for analytic exploration. Then, in 1970s A.A.~Sveshnikov showed that for some special class of SDEs, the right hand side of~\eqref{eq:Pug_gen} can be seen as the image of~$E(\vec{z}, t)$ under the action of a certain singular integral operator, so in fact the equation is integral-differential. Later on, this problem was passed to one of the authors, who found a way to solve it explicitly.

Now, let us see what happens to the general equation~\eqref{eq:Pug_gen} in case of an SDE with piecewise linear coefficients that have two domains of linearity and restate the results from~\cite{Zayats2013}. Assume that the drift vector and the intensity matrix have the form
\begin{equation}
  \label{eq:sde_coef1d_1}
  \vec{a}(t, \vec{x})=\left\{
    \begin{aligned}
      \op{C}^{(1)}(t)\ \vec{x} +  \vec{d}^{(1)}(t), \quad x_1 > 0,\\
      \op{C}^{(2)}(t)\ \vec{x} +  \vec{d}^{(1)}(t), \quad x_1 < 0,
    \end{aligned}
  \right.
\end{equation}
and
\begin{equation}
  \label{eq:sde_coef1d_2}
  \op{H}(t, \vec{x})=\sqrt{2} \left\{
    \begin{aligned}
      \op{H}^{(1)}(t), \quad x_1>0,\\
      \op{H}^{(2)}(t), \quad x_1<0.\\
    \end{aligned}
  \right.
\end{equation}
for some matrices~$\op{C}^{\pm}(t)$,~$\op{H}^{\pm}(t)$ and vectors~$\vec{d}^{\pm}(t)$. Then, the following theorem holds.
\begin{mythm}
  \label{thm1}
  The characteristic function~$E(\vec{z},t)$ of the vector process~$\vec{X}(t)$ given by SDE~\eqref{eq:SDE_gen} with the coefficients~\eqref{eq:sde_coef1d_1} and~\eqref{eq:sde_coef1d_2} satisfies the Pugachev--Sveshnikov equation
  \begin{equation}
    \begin{aligned}
      \label{eq:1pl_pg_gen}
      &\frac{\partial E(\vec{z},t)}{\partial t} + (\vec{z}^{*} \op{B}_0 \vec{z} - i \vec{z}^{*} \vec{d}_0 -  \vec{z}^{*} \op{C}_0 \nabla_{\vec{z}} ) E(\vec{z},t) -\\
      &-i(\vec{z}^{*} \op{B}_1  \vec{z} - i \vec{z}^{*}  \vec{d}_1 -\vec{z}^{*} \op{C}_1 \nabla_{\vec{z}}) \hat{E}(\vec{z},t)=0,
    \end{aligned}
  \end{equation}
  where
  \begin{equation}
    \begin{aligned}
      &\hat{E}(\vec{z},t)=\frac{1}{\pi} \mathrm{v.p.} \int \limits_{- \infty}^{+ \infty} \frac{E|_{z_1=s}}{s-z_1} ds,
    \end{aligned}
  \end{equation}
  and
  \begin{equation}
    \label{eq:1pl_coeff}
    \begin{aligned}
      &\op{C}_0=\frac{\op{C}^{(1)} + \op{C}^{(2)}}{2},\ \op{C}_1=\frac{\op{C}^{(1)} - \op{C}^{(2)}}{2},\\
      &\vec{d}_0=\frac{ \vec{d}^{(1)} +  \vec{d}^{(2)}}{2},\  \vec{d}_1=\frac{ \vec{d}^{(1)} -  \vec{d}^{(2)}}{2},\\
      &\op{H}_0=\frac{\op{H}^{(1)} + \op{H}^{(2)}}{2},\ \op{H}_1=\frac{\op{H}^{(1)} - \op{H}^{(2)}}{2},\\
      & \op{B}_0=\frac{\op{H}^{(1)}(\op{H}^{(1)})^* + \op{H}^{(2)}(\op{H}^{(2)})^*}{2},\\
      &\op{B}_1=\frac{\op{H}^{(1)}(\op{H}^{(1)})^* - \op{H}^{(2)}(\op{H}^{(2)})^*}{2}.
    \end{aligned}
  \end{equation}
\end{mythm}

Not getting into details, it is worth noticing that the main idea of proving~\eqref{eq:1pl_pg_gen} is based on the simple identity
\begin{equation}
  \label{eq:lemma1_fin}
  \mexp{\sign{X}\, e^{i z X}} = \frac{1}{\pi i} \mathrm{v.p.} \int \limits_{-\infty}^{+ \infty} \frac{E(s)}{s-z}\, ds
\end{equation}
for any (absolutely continuous) random variable~$X$ with the characteristic function~$E(z)$.

In order to formulate the next theorem we introduce the Cauchy type integral 
\begin{equation}
  F(\zeta; \vec{z'}, t) = \frac{1}{2 \pi i} \int \limits_{-\infty}^{+\infty} \frac{E|_{z_1=s}}{s - \zeta} ds,\quad \Im{\zeta}\ne 0,
\end{equation}
where~$\vec{z}' = (z_2, \dots, z_n)^* \in \mathbb{R}^{n-1}$,
its limit values as~$\zeta \to z_1 \pm 0 i$, $z_1 \in \mathbb{R}$:
\begin{equation}
  F^\pm(\vec{z}, \tau) = \lim\limits_{\zeta \to z_1 \pm 0 i} F(\zeta, \vec{z'}; t),
\end{equation}
and the following short operator notation:
\begin{equation}
  \begin{aligned}
    &\mathcal{A}^+ = \frac{\partial}{\partial t} + \vec{z}^{*} \op{B}^{(1)}  \vec{z} - i  \vec{z}^{*}  \vec{d}^{(1)} - \vec{z}^{*} \op{C}^{(1)} \nabla_{ \vec{z}},\\
    &\mathcal{A}^- = \frac{\partial}{\partial t} + \vec{z}^{*} \op{B}^{(2)}  \vec{z} - i  \vec{z}^{*}  \vec{d}^{(2)} - \vec{z}^{*} \op{C}^{(2)} \nabla_{ \vec{z}}.
  \end{aligned}
\end{equation}

In above-specified terms the following theorem takes place.
\begin{mythm}
  Assume~$F^\pm(\vec{z},t) = O(1/z_1)$ as~$z_1 \to \infty$ such that~$\Im{z_1} \gtrless 0$. Then the limit values~$F^{\pm}(\vec{z},t)$ satisfy the master equation
  \begin{equation}
    \label{eq:gen_eq_1d}
    \mathcal{A}^{\pm} F^{\pm}(\vec{z},t) = G_0(\vec{z}',t) + G_1(\vec{z}', t) z_1
  \end{equation}
  for some functions~$G_0(\vec{z}',t)$ and $G_1(\vec{z}',t)$.
\end{mythm}

This theorem, in fact, reduces the problem of finding~$F^{\pm}$ to some inverse problem. One need to find the intermediary functions~$G_0$ and~$G_1$, based on analytic properties of the left hand side of~\eqref{eq:gen_eq_1d}, e.g., as it was done in section~\ref{sec:moments}.

Notice, that analogous results can be obtained in a more general setting as well, for instance, when the number of linearity domains is larger~\cite{Zayats2013a}.

\input{ndtcs.bbl}
\end{document}

%% file: ndtcs.bbl
\providecommand{\noopsort}[1]{}\providecommand{\singleletter}[1]{#1}%